\documentclass[12pt]{article}
\usepackage{latexsym,amsfonts,amssymb}
\usepackage[dvips]{color}
\usepackage{colordvi,multicol}
\usepackage{amsmath}
\usepackage{graphicx}
\usepackage{pdflscape}
\usepackage{rotating}

\usepackage{latexsym,amsfonts,amssymb}
\usepackage[dvips]{color}
\usepackage{colordvi,multicol}
\usepackage{amsmath}
\usepackage{graphicx}
\usepackage{pdflscape}
\usepackage{rotating}
\usepackage{breqn}

\newtheorem{thm}{Theorem}[section]

\newtheorem{con}[thm]{Conjecture}

\date{}
\begin{document}

\title{\bf Variation comparison between the $F$-distribution and the normal distribution}
\author{Ping Sun, Ze-Chun Hu  and Wei Sun\thanks{Corresponding author.}\\ \\
  {\small Business School, Chengdu University, Chengdu 610106, China}\\ \\
 {\small College of Mathematics, Sichuan University, Chengdu 610065, China}\\ \\
{\small Department of Mathematics and Statistics, Concordia University, Montreal H3G 1M8,  Canada}\\ \\
{\small  sunping@cdu.edu.cn\ \ \ \ zchu@scu.edu.cn\ \ \ \ wei.sun@concordia.ca}}

\maketitle

\begin{abstract}
\noindent
Let $X_{d_1,d_2}$ be an $F$-random variable with numerator and denominator degrees of freedom $d_1$ and $d_2$, respectively. We investigate the inequality:  $P\{|X_{d_1,d_2}-E[X_{d_1,d_2}]|\le \sqrt{{\rm Var}(X_{d_1,d_2})}\}\ge P\{|W-E[W]|\le \sqrt{{\rm Var}(W)}\}$, where $W$ is a standard normal random variable or a $\chi^2(d_1)$ random variable. We prove that this inequality holds for $d_1\in\{1,2,3,4\}$ and $5\le d_2\in\mathbb{N}$.
\end{abstract}



\section{Introduction}

In \cite{{SHS},SHSb}, we initiated the study of the variation comparison between familiar distributions and the normal distribution. Let $\{X_{\tau}\}$ be a family of random variables with the same type of distributions but with different parameters $\tau$. We discovered that for a large number of infinitely divisible continuous distributions including the Gamma, Log-normal, student's $t$ and inverse Gaussian distributions, we have the following inequality:
$$
P\left\{|X_{\tau}-E[X_{\tau}]|\le \sqrt{{\rm Var}(X_{\tau})}\right\}>P\{|Z|\le 1\}\approx 0.6827,
$$
and
$$
\inf_{\tau}P\left\{|X_{\tau}-E[X_{\tau}]|\le \sqrt{{\rm Var}(X_{\tau})}\right\}=P\{|Z|\le 1\}.
$$
Hereafter $Z$ denotes a standard normal random variable.

To the best of our knowledge, the above variation comparison inequality has not been given in the literature, although many famous distributions have a long history and date back to the era of Bernoulli and De Moivre. We have established in \cite{{SHS},SHSb} the variation comparison inequality for most familiar infinitely divisible continuous distributions listed in Sato's book \cite{S} (cf. \cite[Pages 46 and 47]{S}) except the $F$-distribution. The $F$-distribution is a basic and important probability distribution, which has many statistical applications including comparing the variability of two population samples, comparing the means of two or more groups, and testing the validity of a multiple regression equation.  In this paper, we continue to consider the following conjecture.
\begin{con}\label{con1}
Let $d_1\in \mathbb{N}$, $5\le d_2\in\mathbb{N}$,  $X_{d_1,d_2}$ be an $F$-random variable with parameters $d_1$ and $d_2$, and $Z$ be a standard normal random variable. Then,
$$
P\left\{|X_{d_1,d_2}-E[X_{d_1,d_2}]|\le \sqrt{{\rm Var}(X_{d_1,d_2})}\right\}>P\{|Z|\le 1\}\approx 0.6827,
$$
and
$$
\inf_{d_1,d_2}P\left\{|X_{d_1,d_2}-E[X_{d_1,d_2}]|\le \sqrt{{\rm Var}(X_{d_1,d_2})}\right\}=P\{|Z|\le 1\}.
$$
\end{con}

Different from the two-parameter infinitely divisible distributions considered in \cite{{SHS},SHSb}, it seems not possible to reduce Conjecture \ref{con1} to an inequality problem with just one parameter. This might explain why the variation comparison between the $F$-distribution and the normal distribution is more difficult. We hope the exploration of this paper can inspire us to investigate the related problems for   three-parameter infinitely divisible distributions. We have been able to resolve the conjecture for the case that $d_1$ is small and $d_2$ is arbitrary. More precisely, we have the following result.
\begin{thm}\label{thm25f}
Let $d_1\in \{1,2,3,4\}$, $5\le d_2\in\mathbb{N}$,  $X_{d_1,d_2}$ be an $F$-random variable with parameters $d_1$ and $d_2$,  $Z$ be a standard normal random variable, and $\Gamma_{d_1}$ be a $\chi^2(d_1)$-random variable. Then,
$$
P\left\{|X_{d_1,d_2}-E[X_{d_1,d_2}]|\le \sqrt{{\rm Var}(X_{d_1,d_2})}\right\}>P\{|Z|\le 1\}\approx 0.6827,
$$
and
$$
\inf_{d_2\in\mathbb{N}}P\left\{|X_{d_1,d_2}-E[X_{d_1,d_2}]|\le \sqrt{{\rm Var}(X_{d_1,d_2})}\right\}=P\left\{|\Gamma_{d_1}-E[\Gamma_{d_1}]|\le \sqrt{{\rm Var}(\Gamma_{d_1})}\right\}.
$$
\end{thm}
Note that, by \cite[Theorem 1.2]{SHS},  we have that
$$
\inf_{d_1\in\mathbb{N}}P\left\{|\Gamma_{d_1}-E[\Gamma_{d_1}]|\le \sqrt{{\rm Var}(\Gamma_{d_1})}\right\}=P\{|Z|\le 1\}.
$$

The remainder of this paper is organized as follows. First, we make some preparation in the next section and discuss  different situations that we need to consider for the variation comparison inequality. Then, in Sections 3 and 4, we give the proof of Theorem \ref{thm25f} for cases $d_1=1,2$ and $d_1=3,4$, respectively. Finally, we make some remarks on the case $d_1\ge 5$ in the last section.

\section{Preliminary}\setcounter{equation}{0}

We have (cf. \cite{WF})
$$
E[X_{d_1,d_2}]=\frac{d_2}{d_2-2},\ \ \ \ {\rm Var}(X_{d_1,d_2})=\frac{2d_2^2(d_1+d_2-2)}{d_1(d_2-2)^2(d_2-4)},
$$
and
$$
P\{X_{d_1,d_2}\le x\}=I_{\frac{d_1x}{d_1x+d_2}}\left(\frac{d_1}{2},\frac{d_2}{2}\right),\ \ \ \ x> 0,
$$
where $I$ is the regularized incomplete beta function:
$$
I_x(a,b)=\frac{B(x;a,b)}{B(a,b)}=\frac{\int_0^xt^{a-1}(1-t)^{b-1}dt}{B(a,b)}.
$$
Note that $X_{d_1,d_2}$ has a finite variance only if $d_2>4$. Then, we get
\begin{eqnarray*}
&&P\left\{|X_{d_1,d_2}-E[X_{d_1,d_2}]|\le \sqrt{{\rm Var}(X_{d_1,d_2})}\right\}\\
&=&I_{\frac{d_1}{d_1+(d_2-2)\left(1+\sqrt{\frac{2(d_1+d_2-2)}{d_1(d_2-4)}}\right)^{-1}}}\left(\frac{d_1}{2},\frac{d_2}{2}\right)-I_{\frac{d_1}{d_1+(d_2-2)\left(\max\left\{0,1-\sqrt{\frac{2(d_1+d_2-2)}{d_1(d_2-4)}}\right\}\right)^{-1}}}\left(\frac{d_1}{2},\frac{d_2}{2}\right)\\
&:=&F_{d_1,d_2}.
\end{eqnarray*}
Hereafter, we adopt the convention that $0^{-1}=\infty$ and $\frac{1}{\infty}=0$.

Note that $X_{d_1,d_2}$ has the following expression:
$$
X_{d_1,d_2}=\frac{\Gamma_{d_1}/d_1}{\Gamma_{d_2}/d_2},
$$
where $\Gamma_{d_1}$ and $\Gamma_{d_2}$ are independent $\chi^2$-random variables with respective degrees of freedom $d_1$ and $d_2$. By Slutsky's theorem, we know that
$X_{d_1,d_2}$ converges to $\Gamma_{d_1}/d_1$ in distribution as $d_2\rightarrow\infty$ for any fixed $d_1\in\mathbb{N}$. By \cite[Theorem 1.2]{SHS}, to complete the proof, it suffices to show that
\begin{eqnarray}\label{30A}
F_{d_1,d_2+2}<F_{d_1,d_2},\ \ \ \ \forall d_2\ge5.
\end{eqnarray}

To simplify notation, define
$$
A:={\frac{d_1}{d_1+d_2\left(1+\sqrt{\frac{2(d_1+d_2)}{d_1(d_2-2)}}\right)^{-1}}},\ \ \ \ B:={\frac{d_1}{d_1+(d_2-2)\left(1+\sqrt{\frac{2(d_1+d_2-2)}{d_1(d_2-4)}}\right)^{-1}}},
$$
and
$$
C:={\frac{d_1}{d_1+d_2\left(1-\sqrt{\frac{2(d_1+d_2)}{d_1(d_2-2)}}\right)^{-1}}},\ \ \ \ D:={\frac{d_1}{d_1+(d_2-2)\left(1-\sqrt{\frac{2(d_1+d_2-2)}{d_1(d_2-4)}}\right)^{-1}}}.
$$
By the  relation $I_x(a,b+1)=I_x(a,b)+\frac{x^a(1-x)^b}{bB(a,b)}$, we get
\begin{eqnarray}\label{May1311}
&&(\ref{30A})\ {\rm holds}\nonumber\\
&\Leftrightarrow&\ \ I_A\left(\frac{d_1}{2},\frac{d_2}{2}\right)+\frac{A^{\frac{d_1}{2}}(1-A)^{\frac{d_2}{2}}}{\frac{d_2}{2}B\left(\frac{d_1}{2},\frac{d_2}{2}\right)}- I_C\left(\frac{d_1}{2},\frac{d_2}{2}\right)-\frac{C^{\frac{d_1}{2}}(1-C)^{\frac{d_2}{2}}}{\frac{d_2}{2}B\left(\frac{d_1}{2},\frac{d_2}{2}\right)}\nonumber\\
&&<I_B\left(\frac{d_1}{2},\frac{d_2}{2}\right)-I_D\left(\frac{d_1}{2},\frac{d_2}{2}\right)\nonumber\\
&\Leftrightarrow&\ \ 2A^{\frac{d_1}{2}}(1-A)^{\frac{d_2}{2}}+d_2\int_C^Dt^{\frac{d_1}{2}-1}(1-t)^{\frac{d_2}{2}-1}dt\nonumber\\
&&<d_2\int_A^Bt^{\frac{d_1}{2}-1}(1-t)^{\frac{d_2}{2}-1}dt+2C^{\frac{d_1}{2}}(1-C)^{\frac{d_2}{2}}.
\end{eqnarray}

Note that $A<B$
since
$$
{\frac{d_1}{d_1+d_2\left(1+\sqrt{\frac{2(d_1+d_2)}{d_1(d_2-2)}}\right)^{-1}}}
$$
is a strictly decreasing function of $d_2$ for any fixed $d_1$. We have
\begin{eqnarray}
&&C>0\Leftrightarrow1-\sqrt{\frac{2(d_1+d_2)}{d_1(d_2-2)}}>0,\nonumber\\
&&D>0\Leftrightarrow1-\sqrt{\frac{2(d_1+d_2-2)}{d_1(d_2-4)}}>0,\label{34}\\
&&D>0\Rightarrow C>0.\nonumber
\end{eqnarray}
We would like to point out that  $C>D$ is possible. Under the assumption that $D>0$, we have
\begin{eqnarray}\label{35}
D>C&\Leftrightarrow&d_2\left\{1-\sqrt{\frac{2(d_1+d_2)}{d_1(d_2-2)}}\right\}^{-1}>(d_2-2)\left\{1-\sqrt{\frac{2(d_1+d_2-2)}{d_1(d_2-4)}}\right\}^{-1}\nonumber\\
&\Leftrightarrow&d_2\left\{1-\sqrt{\frac{2(d_1+d_2-2)}{d_1(d_2-4)}}\right\}>(d_2-2)\left\{1-\sqrt{\frac{2(d_1+d_2)}{d_1(d_2-2)}}\right\}\nonumber\\
&\Leftrightarrow&2+(d_2-2)\sqrt{\frac{2(d_1+d_2)}{d_1(d_2-2)}}>d_2\sqrt{\frac{2(d_1+d_2-2)}{d_1(d_2-4)}}\nonumber\\
&\Leftrightarrow&4+\frac{2(d_1+d_2)(d_2-2)}{d_1}+4(d_2-2)\sqrt{\frac{2(d_1+d_2)}{d_1(d_2-2)}}>\frac{2d_2^2(d_1+d_2-2)}{d_1(d_2-4)}\nonumber\\
&\Leftrightarrow&4(d_2-2)\sqrt{\frac{2(d_1+d_2)}{d_1(d_2-2)}}>\frac{8d_2(d_1+d_2-2)}{d_1(d_2-4)}\nonumber\\
&\Leftrightarrow&(d_2-4)\sqrt{d_1(d_2-2)(d_1+d_2)}>\sqrt{2}d_2(d_1+d_2-2)\nonumber\\
&\Leftrightarrow&d_1(d_2-2)(d_1+d_2)(d_2-4)^2>2d^2_2(d_1+d_2-2)^2.
\end{eqnarray}
If $d_1\ge 3$, then we can show that  actually (\ref{35}) holds without the extra assumption that $D>0$.  That is, we have that
\begin{eqnarray}\label{78}
D>C\Leftrightarrow d_1(d_2-2)(d_1+d_2)(d_2-4)^2>2d^2_2(d_1+d_2-2)^2,\ \ \ \ \forall d_1\ge3.
\end{eqnarray}
In fact, suppose the right hand side of (\ref{78}) holds. By  (\ref{34}), we get
\begin{eqnarray*}
D>0&\Leftrightarrow& {d_1^2d_2^2(d_2-4)^2}>4d_2^2(d_1+d_2-2)^2\\
&\Leftarrow&{d_1^2d_2^2(d_2-4)^2}>2d_1(d_2-2)(d_1+d_2)(d_2-4)^2\\
&\Leftrightarrow&d_1d_2^2>2(d_2-2)(d_1+d_2)\\
&\Leftrightarrow&(d_1-2)d^2_2+4d_1+4d_2>2d_1d_2\\
&\Leftarrow&(d_1-2)d_2+4>2d_1\\
&\Leftarrow&(d_1-2)5+4>2d_1\\
&\Leftrightarrow&d_1>2.
\end{eqnarray*}

There are three different situations when considering  inequality (\ref{May1311}).

\noindent Condition 1:
\begin{eqnarray*}
C=0\Leftrightarrow1-\sqrt{\frac{2(d_1+d_2)}{d_1(d_2-2)}}\le0\Leftrightarrow\left\{ \begin{array}{ll}
        d_2\ge 5,\ \  \ \ & \mbox{if $d_1=1$ or 2},\\
        5\le d_2\le 4+\frac{8}{d_1-2},\ \ \ \  & \mbox{if $3\le d_1\le 10$};\end{array} \right.
\end{eqnarray*}
\noindent Condition 2:
\begin{eqnarray*}
&&C>0\ {\rm and }\ D= 0\\
&\Leftrightarrow&1-\sqrt{\frac{2(d_1+d_2)}{d_1(d_2-2)}}>0\ {\rm and}\ 1-\sqrt{\frac{2(d_1+d_2-2)}{d_1(d_2-4)}}\le 0\\
&\Leftrightarrow&\left\{ \begin{array}{ll}
        4+\frac{8}{d_1-2}< d_2\le6+\frac{8}{d_1-2},\ \  \ \ & \mbox{if $3\le d_1\le 10$},\\
        d_2=5\ {\rm or}\ 6,\ \ \ \  & \mbox{if $d_1\ge 11$};\end{array} \right.
\end{eqnarray*}
Condition 3:
\begin{eqnarray}\label{cond3}
D>0&\Leftrightarrow&1-\sqrt{\frac{2(d_1+d_2-2)}{d_1(d_2-4)}}\ge0\Leftrightarrow d_1\ge 3\ {\rm and}\ d_2> 6+\frac{8}{d_1-2}.
\end{eqnarray}

\section{Proof of  (\ref{30A}) for cases $d_1=1,2$}\setcounter{equation}{0}

For these two cases, $C=D=0$. By (\ref{May1311}), it suffices to prove the following inequality:
\begin{eqnarray}\label{May13aa}
2A^{\frac{d_1}{2}}(1-A)^{\frac{d_2}{2}}<d_2\int_A^Bt^{\frac{d_1}{2}-1}(1-t)^{\frac{d_2}{2}-1}dt.
\end{eqnarray}
By virtue of the transformation $t=Aw$, we get
\begin{eqnarray}\label{L1}
&&(\ref{May13aa})\ {\rm holds}\nonumber\\
&\Leftrightarrow&2(1-A)<d_2\int_1^{\frac{B}{A}}w^{\frac{d_1}{2}-1}\left(1-\frac{A}{1-A}(w-1)\right)^{\frac{d_2}{2}-1}dw\nonumber\\
&\Leftrightarrow&2(1-A)<d_2\int_0^{\frac{B-A}{A}}(1+y)^{\frac{d_1}{2}-1}\left(1-\frac{A}{1-A}y\right)^{\frac{d_2}{2}-1}dy\nonumber\\
&\Leftrightarrow&2A(1-A)<d_2(B-A)\int_0^1\left(1+\frac{B-A}{A}x\right)^{\frac{d_1}{2}-1}\left(1-\frac{B-A}{1-A}x\right)^{\frac{d_2}{2}-1}dx.
\end{eqnarray}

\subsection{Case $d_1=2$}

Suppose  $d_1=2$. By (\ref{L1}), we get
\begin{eqnarray*}
(\ref{May13aa})\ {\rm holds}
&\Leftrightarrow&2A(1-A)<2(1-A)\left[1-\left(\frac{1-B}{1-A}\right)^{\frac{d_2}{2}}\right]\nonumber\\
&\Leftrightarrow&(1-B)^{\frac{d_2}{2}}<(1-A)^{\frac{d_2}{2}+1}.
\end{eqnarray*}
We have
\begin{eqnarray}\label{L11}
&&(1-B)^{\frac{d_2}{2}}<(1-A)^{\frac{d_2}{2}+1}\nonumber\\
&\Leftrightarrow&\left[{\frac{d_2-2}{d_1\left(1+\sqrt{\frac{2(d_1+d_2-2)}{d_1(d_2-4)}}\right)+(d_2-2)}}\right]^{\frac{d_2}{2}}<\left[{\frac{d_2}{d_1\left(1+\sqrt{\frac{2(d_1+d_2)}{d_1(d_2-2)}}\right)+d_2}}\right]^{\frac{d_2}{2}+1}\nonumber\\
&\Leftrightarrow&{\frac{d_2}{2}}\ln\left[{\frac{d_2-2}{d_1\left(1+\sqrt{\frac{2(d_1+d_2-2)}{d_1(d_2-4)}}\right)+(d_2-2)}}\right]<\left({\frac{d_2}{2}+1}\right)\ln\left[{\frac{d_2}{d_1\left(1+\sqrt{\frac{2(d_1+d_2)}{d_1(d_2-2)}}\right)+d_2}}\right]\nonumber\\
&\Leftrightarrow&{\frac{d_2}{2}}\ln\left[1+\frac{d_1}{d_2-2}\left(1+\sqrt{\frac{2(d_1+d_2-2)}{d_1(d_2-4)}}\right)\right]>\left({\frac{d_2}{2}+1}\right)\ln\left[1+\frac{d_1}{d_2}\left(1+\sqrt{\frac{2(d_1+d_2)}{d_1(d_2-2)}}\right)\right].\nonumber\\
&&
\end{eqnarray}

For $x=1,2,3$ and $3\le y\in\mathbb{R}$, define
\begin{eqnarray}\label{L12}
H_x(y):=\left({\frac{y}{2}+1}\right)\ln\left[1+\frac{x}{y}\left(1+\sqrt{\frac{2(x+y)}{x(y-2)}}\right)\right].
\end{eqnarray}
We will prove  (\ref{L11}) by showing that $H_2$ is a decreasing function for $y\in\{3,4,5\}$ and
$H'_2(y)<0$ for $5\le y\in\mathbb{R}$. By virtue of Mathematica, we obtain  the following values, which are decreasing.
\begin{eqnarray*}
H_2(3)=2.87436,\ \ H_2(4)=2.58363,\ \ H_2(5)=2.44523.
\end{eqnarray*}

We have
\begin{eqnarray}\label{L13}
H'_x(y)=\frac{1}{2}\ln\left[1+\frac{x}{y}\left(1+\sqrt{\frac{2(x+y)}{x(y-2)}}\right)\right]+\left({\frac{y}{2}+1}\right)\left\{\ln\left[1+\frac{x}{y}\left(1+\sqrt{\frac{2(x+y)}{x(y-2)}}\right)\right]\right\}'.
\end{eqnarray}
Now let $x=2$. Define
$$
p_2(y):=\frac{2}{y}\left(1+\sqrt{\frac{2+y}{y-2}}\right).
$$
We have
\begin{eqnarray*}
H'_2(y)
&\le&\frac{1}{2}\left(p_2(y)-\frac{p_2^2(y)}{2}+\frac{p_2^3(y)}{3}\right)+\left({\frac{y}{2}+1}\right)\left\{\ln\left[1+\frac{2}{y}\left(1+\sqrt{\frac{2+y}{y-2}}\right)\right]\right\}'\\
&:=&L_2(y).
\end{eqnarray*}
By virtue of Mathematica, we obtain that for $y\ge 5$,
\begin{eqnarray*}
&&-3y^3(2+y) \left(2+y+2\sqrt{\frac{2+y}{y-2}}\right)L_2(y)\\
&=&2\left(\frac{2+y}{y-2}\right)^{\frac{3}{2}}\left[32\sqrt{\frac{2+y}{y-2}}-4y\left(9+5\sqrt{\frac{2+y}{y-2}}\right)+2y^2\left(2-\sqrt{\frac{2+y}{y-2}}\right)+3y^3\right]\\
&>&2\left(\frac{2+y}{y-2}\right)^{\frac{3}{2}}y\left(39-20\sqrt{\frac{2+y}{y-2}}\right)\\
&=&\frac{1}{20}\left(\frac{2+y}{y-2}\right)^{\frac{3}{2}}y\left(1.95-\sqrt{1+\frac{4}{y-2}}\right)\\
&\ge&\frac{1}{20}\left(\frac{2+y}{y-2}\right)^{\frac{3}{2}}y\left(1.95-\sqrt{1+\frac{4}{5-2}}\right)\\
&>& 0.42247\cdot\frac{1}{20}\left(\frac{2+y}{y-2}\right)^{\frac{3}{2}}y,
\end{eqnarray*}
Then, $H'_2(y)<0$ for $5\le y\in\mathbb{R}$. Thus, $(1-B)^{\frac{d_2}{2}}<(1-A)^{\frac{d_2}{2}+1}$ by (\ref{L11}). Therefore, (\ref{May13aa}) holds for the case $d_1=2$.

\subsection{Case $d_1=1$}

Suppose $d_1=1$. By (\ref{L1}), we get
\begin{eqnarray*}
&&(\ref{May13aa})\ {\rm holds}\nonumber\\
&\Leftarrow&2A(1-A)<d_2(B-A)\int_0^1\left(1+\frac{-\frac{1}{2}(B-A)}{A}x\right)\left(1-\frac{B-A}{1-A}x\right)^{\frac{d_2}{2}-1}dx\nonumber\\
&\Leftrightarrow&2A(1-A)<d_2(B-A)\left\{\frac{3A-1}{2A}\int_0^1\left(1-\frac{B-A}{1-A}x\right)^{\frac{d_2}{2}-1}dx+\frac{1-A}{2A}\int_0^1\left(1-\frac{B-A}{1-A}x\right)^{\frac{d_2}{2}}dx\right\}
\end{eqnarray*}

\begin{eqnarray}\label{gfd}
&\Leftrightarrow&2A(1-A)<\frac{(3A-1)(1-A)}{A}\left[1-\left(\frac{1-B}{1-A}\right)^{\frac{d_2}{2}}\right]+\frac{d_2(1-A)^2}{(d_2+2)A}\left[1-\left(\frac{1-B}{1-A}\right)^{\frac{d_2}{2}+1}\right]\nonumber\\
&\Leftrightarrow&2A^2<\frac{2(d_2+3)A-2}{d_2+2}-\frac{3(d_2+2)A-2-d_2B}{d_2+2}\left(\frac{1-B}{1-A}\right)^{\frac{d_2}{2}}\nonumber\\
&\Leftrightarrow&\left[3(d_2+2)A-2-d_2B\right](1-B)^{\frac{d_2}{2}}<2[(d_2+2)A-1](1-A)^{\frac{d_2}{2}+1}.
\end{eqnarray}

Below we analyze the coefficients of inequality (\ref{gfd}). First, we  have
\begin{eqnarray*}
(d_2+2)A>1\Leftrightarrow d_1(d_2+2)>{d_1+d_2\left(1+\sqrt{\frac{2(d_1+d_2)}{d_1(d_2-2)}}\right)^{-1}},
\end{eqnarray*}
which holds obviously. Second, we have
\begin{eqnarray*}
&&3(d_2+2)A-2-d_2B>0\\
&\Leftrightarrow&3(d_2+2)A>2+d_2B\\
&\Leftrightarrow&\frac{3d_1(d_2+2)\left(1+\sqrt{\frac{2(d_1+d_2)}{d_1(d_2-2)}}\right)}{d_1\left(1+\sqrt{\frac{2(d_1+d_2)}{d_1(d_2-2)}}\right)+d_2}>{\frac{d_1(d_2+2)\left(1+\sqrt{\frac{2(d_1+d_2-2)}{d_1(d_2-4)}}\right)+2(d_2-2)}{d_1\left(1+\sqrt{\frac{2(d_1+d_2-2)}{d_1(d_2-4)}}\right)+(d_2-2)}}\\
&\Leftrightarrow&\ \ 2d_1^2(d_2+2)\left(1+\sqrt{\frac{2(d_1+d_2)}{d_1(d_2-2)}}\right)\left(1+\sqrt{\frac{2(d_1+d_2-2)}{d_1(d_2-4)}}\right)\\
&&\ \ +d_1(d_2-2)(3d_2+4)\left(1+\sqrt{\frac{2(d_1+d_2)}{d_1(d_2-2)}}\right)\\
&&>d_1d_2(d_2+2)\left(1+\sqrt{\frac{2(d_1+d_2-2)}{d_1(d_2-4)}}\right)+2d_2(d_2-2)\\
&\Leftrightarrow&\ \ 2(d_2-2)+(3d^2_2-4)\sqrt{\frac{2(1+d_2)}{d_2-2}}+2(d_2+2)\left(1+\sqrt{\frac{2(1+d_2)}{d_2-2}}\right)\sqrt{\frac{2(d_2-1)}{d_2-4}}\\
&&>d_2(d_2+2)\sqrt{\frac{2(d_2-1)}{d_2-4}}\\
&\Leftarrow&d_2(3d_2+2)\sqrt{\frac{2(1+d_2)}{d_2-2}}>d_2(d_2+2)\sqrt{\frac{2(d_2-1)}{d_2-4}}\\
&\Leftrightarrow&(3d_2+2)^2(1+d_2)(d_2-4)>(d_2+2)^2(d_2-1)(d_2-2)\\
&\Leftarrow&2.4^2(d_2-4)>d_2-2\\
&\Leftarrow&5(d_2-4)>d_2-2\\
&\Leftrightarrow&4d_2>18,
\end{eqnarray*}
which holds for $d_2\ge 5$. Third, we have
\begin{eqnarray}\label{ret}
&&2[(d_2+2)A-1]>3(d_2+2)A-2-d_2B\nonumber\\
&\Leftrightarrow&d_2B>(d_2+2)A\nonumber\\
&\Leftrightarrow&{\frac{d_1d_2\left(1+\sqrt{\frac{2(d_1+d_2-2)}{d_1(d_2-4)}}\right)}{d_1\left(1+\sqrt{\frac{2(d_1+d_2-2)}{d_1(d_2-4)}}\right)+(d_2-2)}}>\frac{d_1(d_2+2)\left(1+\sqrt{\frac{2(d_1+d_2)}{d_1(d_2-2)}}\right)}{d_1\left(1+\sqrt{\frac{2(d_1+d_2)}{d_1(d_2-2)}}\right)+d_2}\nonumber\\
&\Leftrightarrow&{\frac{d_2\left(1+\sqrt{\frac{2(d_2-1)}{d_2-4}}\right)}{\left(1+\sqrt{\frac{2(d_2-1)}{d_2-4}}\right)+(d_2-2)}}>\frac{(d_2+2)\left(1+\sqrt{\frac{2(1+d_2)}{d_2-2}}\right)}{\left(1+\sqrt{\frac{2(1+d_2)}{d_2-2}}\right)+d_2}.
\end{eqnarray}
Define
$$
K(y):={\frac{y\left(1+\sqrt{\frac{2(y-1)}{y-4}}\right)}{\left(1+\sqrt{\frac{2(y-1)}{y-4}}\right)+(y-2)}},\ \ \ \ 5\le y\in\mathbb{R}.
$$
By virtue of Mathematica, we get,
\begin{eqnarray*}
&&-2(y-4)^2\sqrt{\frac{y-1}{y-4}}\left(y-1+\sqrt{\frac{2(y-1)}{y-4}}\right)^2K'(y)\\
&=&\left(3\sqrt{2}-2\sqrt{\frac{y-1}{y-4}}\right)y^2-\left(6\sqrt{2}-4\sqrt{\frac{y-1}{y-4}}\right)y+16\sqrt{\frac{y-1}{y-4}}\\
&>&\left(9\sqrt{2}-6\sqrt{\frac{y-1}{y-4}}\right)y,
\end{eqnarray*}
which is positive for $y\ge 5$. Then, $K'(y)<0$ for $y\ge 5$. Thus, by (\ref{ret}), we obtain that
$$
2[(d_2+2)A-1]>3(d_2+2)A-2-d_2B.
$$
Hence, by (\ref{gfd}), we conclude that
\begin{eqnarray*}
(\ref{May13aa})\ {\rm holds}\Leftarrow(1-B)^{\frac{d_2}{2}}<(1-A)^{\frac{d_2}{2}+1}.
\end{eqnarray*}

Let $H_x(y)$ be defined by  (\ref{L12}). Note that  (\ref{L11}) and  (\ref{L13}) hold. Now let $x=1$. We have
\begin{eqnarray*}
H'_1(y)
&=&\frac{1}{2}\ln\left[1+\frac{1}{y}\left(1+\sqrt{\frac{2(1+y)}{y-2}}\right)\right]+\left({\frac{y}{2}+1}\right)\left\{\ln\left[1+\frac{1}{y}\left(1+\sqrt{\frac{2(1+y)}{y-2}}\right)\right]\right\}'\\
&:=&L_1(y).
\end{eqnarray*}
Define
$$
p_1(y):=\frac{1}{y}\left(1+\sqrt{\frac{2(1+y)}{y-2}}\right).
$$
We have
\begin{eqnarray*}
H'_1(y)
&\le&\frac{1}{2}\left(p_1(y)-\frac{p_1^2(y)}{2}+\frac{p_1^3(y)}{3}\right)+\left({\frac{y}{2}+1}\right)\left\{\ln\left[1+\frac{1}{y}\left(1+\sqrt{\frac{2(1+y)}{y-2}}\right)\right]\right\}'\\
&:=&L_1(y).
\end{eqnarray*}
By virtue of Mathematica, we obtain that for $y\ge 3$,
\begin{eqnarray*}
&&-12 y^3(y-2)^2 \sqrt{\frac{1+y}{y-2}}\left(1 +y+\sqrt{\frac{2(1+y)}{y-2}} \right)L_1(y)\\
&=&32\sqrt{\frac{1+y}{y-2}}-4y\left(7\sqrt{2}-2\sqrt{\frac{1+y}{y-2}}\right)-y^2\left(35\sqrt{2}-4\sqrt{\frac{1+y}{y-2}}\right)\\
&&+y^3\left(17\sqrt{2}-23\sqrt{\frac{1+y}{y-2}}\right)+3y^4\left(5\sqrt{2}+\sqrt{\frac{1+y}{y-2}}\right)\\
&>&\sqrt{2}y(-28-35y+16y^2)+y^3\left(46\sqrt{2}-23\sqrt{\frac{1+y}{y-2}}\right)\\
&>&0.
\end{eqnarray*}
Then, $H'_1(y)<0$ for $3\le y\in\mathbb{R}$. Thus, $(1-B)^{\frac{d_2}{2}}<(1-A)^{\frac{d_2}{2}+1}$ by (\ref{L11}). Therefore, (\ref{May13aa}) holds for the case $d_1=1$.

\section{Proof of  (\ref{30A}) for cases $d_1=3,4$}\setcounter{equation}{0}

By (\ref{May1311}), it suffices to prove the following two inequalities:
\begin{eqnarray}
&&2A^{\frac{d_1}{2}}(1-A)^{\frac{d_2}{2}}<d_2\int_A^Bt^{\frac{d_1}{2}-1}(1-t)^{\frac{d_2}{2}-1}dt,\label{May14z1}\\
&&2C^{\frac{d_1}{2}}(1-C)^{\frac{d_2}{2}}>d_2\int_C^Dt^{\frac{d_1}{2}-1}(1-t)^{\frac{d_2}{2}-1}dt.\label{May14z2}
\end{eqnarray}

\subsection{Case $d_1=4$}

\subsubsection{Proof of  (\ref{May14z1})}

By (\ref{L1}), we get
\begin{eqnarray*}
&&(\ref{May14z1})\ {\rm holds}\nonumber\\
&\Leftrightarrow&2A(1-A)<d_2(B-A)\int_0^1\left(1+\frac{B-A}{A}x\right)\left(1-\frac{B-A}{1-A}x\right)^{\frac{d_2}{2}-1}dx\nonumber\\
&\Leftrightarrow&2A(1-A)<d_2(B-A)\left\{\frac{1}{A}\int_0^1\left(1-\frac{B-A}{1-A}x\right)^{\frac{d_2}{2}-1}dx-\frac{1-A}{A}\int_0^1\left(1-\frac{B-A}{1-A}x\right)^{\frac{d_2}{2}}dx\right\}
\end{eqnarray*}

\begin{eqnarray}\label{RTD}
&\Leftrightarrow&2A(1-A)<\left\{\frac{2(1-A)}{A}\left[1-\left(\frac{1-B}{1-A}\right)^{\frac{d_2}{2}}\right]-\frac{2d_2(1-A)^2}{(d_2+2)A}\left[1-\left(\frac{1-B}{1-A}\right)^{\frac{d_2}{2}+1}\right]\right\}\nonumber\\
&\Leftrightarrow&(d_2+2)A^2<\left\{(d_2A+2)-(d_2B+2)\left(\frac{1-B}{1-A}\right)^{\frac{d_2}{2}}\right\}\nonumber\\
&\Leftrightarrow&(d_2+2)A^2(1-A)^{\frac{d_2}{2}}<\left\{(d_2A+2)(1-A)^{\frac{d_2}{2}}-(d_2B+2)(1-B)^{\frac{d_2}{2}}\right\}\nonumber\\
&\Leftrightarrow&(d_2B+2)(1-B)^{\frac{d_2}{2}}<[(d_2+2)A+2](1-A)^{\frac{d_2}{2}+1}\nonumber\\
&\Leftrightarrow&\ln(d_2B+2)+\frac{d_2}{2}\ln(1-B)<\ln[(d_2+2)A+2]+\left({\frac{d_2}{2}+1}\right)\ln (1-A)\nonumber\\
&\Leftrightarrow&\ \ \ln\left[\frac{4(d_2+2)\left(1+\sqrt{\frac{2+d_2}{2(d_2-4)}}\right)+2(d_2-2)}{4\left(1+\sqrt{\frac{2+d_2}{2(d_2-4)}}\right)+(d_2-2)}\right]-\frac{d_2}{2}\ln\left[1+\frac{4}{d_2-2}\left(1+\sqrt{\frac{2+d_2}{2(d_2-4)}}\right)\right]\nonumber\\
&&<\ln\left[\frac{4(d_2+4)\left(1+\sqrt{\frac{4+d_2}{2(d_2-2)}}\right)+2d_2}{4\left(1+\sqrt{\frac{4+d_2}{2(d_2-2)}}\right)+d_2}\right]-\frac{d_2+2}{2}\ln\left[1+\frac{4}{d_2}\left(1+\sqrt{\frac{4+d_2}{2(d_2-2)}}\right)\right].
\end{eqnarray}

For $3\le y\in\mathbb{R}$, define
\begin{eqnarray*}
H_4(y)&:=&-\frac{y+2}{2}\ln\left[1+\frac{4}{y}\left(1+\sqrt{\frac{4+y}{2(y-2)}}\right)\right]+\ln\left[\frac{4(y+4)\left(1+\sqrt{\frac{4+y}{2(y-2)}}\right)+2y}{4\left(1+\sqrt{\frac{4+y}{2(y-2)}}\right)+y}\right]\\
&=&\frac{y+2}{2}\ln y-\frac{y+4}{2}\ln\left[4\left(1+\sqrt{\frac{4+y}{2(y-2)}}\right)+y\right]+\ln\left[{4(y+4)\left(1+\sqrt{\frac{4+y}{2(y-2)}}\right)+2y}\right].
\end{eqnarray*}
We will prove  (\ref{RTD}) by showing that $H_4$ is an increasing function for $y\in\{3,4,\dots,12\}$ and
$H'_4(y)>0$ for $12\le y\in\mathbb{R}$. By virtue of Mathematica, we obtain  the following values, which are increasing.
\begin{eqnarray*}
&&H_4(3)=-2.15017,\ \ H_4(4)=-1.85244,\ \ H_4(5)=-1.70932,\ \ H_4(6)=-1.62225,\\
&&H_4(7)=-1.563,\ \ \ \ H_4(8)=-1.51983,\ \ H_4(9)=-1.48688,\ \ H_4(10)=-1.46088,\\
&&H_4(11)=-1.43981,\ \ H_4(12)=-1.42238.
\end{eqnarray*}

We have
\begin{eqnarray*}
H'_4(y)&=&\frac{1}{2}\ln y+\frac{y+2}{2y}\\
&&-\frac{1}{2}\ln\left[4\left(1+\sqrt{\frac{4+y}{2(y-2)}}\right)+y\right]-\frac{y+4}{2}\left\{\ln\left[4\left(1+\sqrt{\frac{4+y}{2(y-2)}}\right)+y\right]\right\}'\\
&&+\left\{\ln\left[{4(y+4)\left(1+\sqrt{\frac{4+y}{2(y-2)}}\right)+2y}\right]\right\}'\\
&=&-\frac{1}{2}\ln\left[1+\frac{4}{y}\left(1+\sqrt{\frac{4+y}{2(y-2)}}\right)\right]+\frac{y+2}{2y}\\
&&-\frac{y+4}{2}\left\{\ln\left[4\left(1+\sqrt{\frac{4+y}{2(y-2)}}\right)+y\right]\right\}'+\left\{\ln\left[{4(y+4)\left(1+\sqrt{\frac{4+y}{2(y-2)}}\right)+2y}\right]\right\}'.
\end{eqnarray*}
Define
$$
p_4(y):=\frac{4}{y}\left(1+\sqrt{\frac{4+y}{2(y-2)}}\right).
$$
We have
\begin{eqnarray*}
-H'_4(y)
&\le&\frac{1}{2}\left(p_4(y)-\frac{p_4^2(y)}{2}+\frac{p_4^3(y)}{3}\right)-\frac{y+2}{2y}\\
&&+\frac{y+4}{2}\left\{\ln\left[4\left(1+\sqrt{\frac{4+y}{2(y-2)}}\right)+y\right]\right\}'-\left\{\ln\left[{4(y+4)\left(1+\sqrt{\frac{4+y}{2(y-2)}}\right)+2y}\right]\right\}'\\
&:=&L_4(y).
\end{eqnarray*}
By virtue of Mathematica, we get
\begin{eqnarray*}
&&-3y^3 (y-2)^2  \left(4+y+2\sqrt{\frac{2(4+y)}{y-2}}\right)\left[8+4\sqrt{\frac{2(4+y)}{y-2}}+y\left(3+\sqrt{\frac{2(4+y)}{y-2}}\right)\right]L_4(y)\\
&=&4\sqrt{\frac{4+y}{y-2}}\Bigg[(2048 + 1536 y - 1040 y^2 - 368 y^3 + 10 y^4 + 3 y^5)\sqrt{2}\\
&&\ \ \ \ \ \ \ \ \ \ \ \ \ \ \ +(-2048 + 3072 y - 112 y^2 - 484 y^3 + 8 y^4 + 3 y^5)\sqrt{\frac{4+y}{y-2}}\Bigg].
\end{eqnarray*}

Define
\begin{eqnarray*}
T_{1}(y)&:=&2048 + 1536 y - 1040 y^2 - 368 y^3 + 10 y^4 + 3 y^5,\\
T_{2}(y)&:=&-2048 + 3072 y - 112 y^2 - 484 y^3 + 8 y^4 + 3 y^5.
\end{eqnarray*}
We have
\begin{eqnarray*}
T_{1}(12+r)&=&188672 + 197760 r + 46192 r^2 + 4432 r^3 + 190 r^4 + 3 r^5,\\
T_{2}(12+r)&=&94720 + 157632 r + 41216 r^2 + 4220 r^3 + 188 r^4 + 3 r^5.
\end{eqnarray*}
Then, $T_{1}(y),T_{2}(y)>0$ for $12\le y\in\mathbb{R}$. Thus, $H'_4(y)>0$ for $12\le y\in\mathbb{R}$. Therefore, (\ref{RTD}) holds for the case  $d_1=4$.

\subsubsection{Proof of  (\ref{May14z2})}

We assume without loss of generality that $D>C>0$. By (\ref{cond3}), we get $D>0\Leftrightarrow d_2\ge 11$. Further, by (\ref{78}), we get
\begin{eqnarray*}
D>C>0&\Leftrightarrow&2(d_2-2)(d_2+4)(d_2-4)^2>d^2_2(d_2+2)^2\\
&\Leftrightarrow&-256 + 192 d_2 - 20 d_2^2 - 16 d_2^3 + d_2^4>0.
\end{eqnarray*}
Define
$$
P_4(d_2):=-256 + 192 d_2 - 20 d_2^2 - 16 d_2^3 + d_2^4.
$$
We have
\begin{eqnarray*}
&&P_4(11)=-7219,\ \ P_4(12)=-7744,\ \ P_4(13)=-7731,\\
&&P_4(14)=-6976,\ \ P_4(15)=-5251,\ \ P_4(16)=-2304,
\\
&&P_4(17+r)=2141 + 5292 r + 898 r^2 + 52 r^3 + r^4.
\end{eqnarray*}
Then,
$$
D>C>0\Leftrightarrow d_2\ge 17.
$$

Similar to (\ref{RTD}), we get
\begin{eqnarray}\label{RTD2}
&&(\ref{May14z2})\ {\rm holds}\nonumber\\
&\Leftrightarrow&\ln(d_2D+2)+\frac{d_2}{2}\ln(1-D)>\ln[(d_2+2)C+2]+\left({\frac{d_2}{2}+1}\right)\ln (1-C)\nonumber\\
&\Leftrightarrow&\ \ \ln\left[\frac{4(d_2+2)\left(1-\sqrt{\frac{2+d_2}{2(d_2-4)}}\right)+2(d_2-2)}{4\left(1-\sqrt{\frac{2+d_2}{2(d_2-4)}}\right)+(d_2-2)}\right]-\frac{d_2}{2}\ln\left[1+\frac{4}{d_2-2}\left(1-\sqrt{\frac{2+d_2}{2(d_2-4)}}\right)\right]\nonumber\\
&&>\ln\left[\frac{4(d_2+4)\left(1-\sqrt{\frac{4+d_2}{2(d_2-2)}}\right)+2d_2}{4\left(1-\sqrt{\frac{4+d_2}{2(d_2-2)}}\right)+d_2}\right]-\frac{d_2+2}{2}\ln\left[1+\frac{4}{d_2}\left(1-\sqrt{\frac{4+d_2}{2(d_2-2)}}\right)\right].\ \ \ \ \ \
\end{eqnarray}
For $15\le y\in\mathbb{R}$, define
\begin{eqnarray*}
R_4(y)&:=&-\frac{y+2}{2}\ln\left[1+\frac{4}{y}\left(1-\sqrt{\frac{4+y}{2(y-2)}}\right)\right]+\ln\left[\frac{4(y+4)\left(1-\sqrt{\frac{4+y}{2(y-2)}}\right)+2y}{4\left(1-\sqrt{\frac{4+y}{2(y-2)}}\right)+y}\right]\\
&=&\frac{y+2}{2}\ln y-\frac{y+4}{2}\ln\left[4\left(1-\sqrt{\frac{4+y}{2(y-2)}}\right)+y\right]+\ln\left[{4(y+4)\left(1-\sqrt{\frac{4+y}{2(y-2)}}\right)+2y}\right].
\end{eqnarray*}
We will prove  (\ref{RTD2}) by showing that
$R'_4(y)<0$ for $15\le y\in\mathbb{R}$.

We have
\begin{eqnarray*}
R'_4(y)&=&\frac{1}{2}\ln y+\frac{y+2}{2y}\\
&&-\frac{1}{2}\ln\left[4\left(1-\sqrt{\frac{4+y}{2(y-2)}}\right)+y\right]-\frac{y+4}{2}\left\{\ln\left[4\left(1-\sqrt{\frac{4+y}{2(y-2)}}\right)+y\right]\right\}'\\
&&+\left\{\ln\left[{4(y+4)\left(1-\sqrt{\frac{4+y}{2(y-2)}}\right)+2y}\right]\right\}'\\
&=&-\frac{1}{2}\ln\left[1+\frac{4}{y}\left(1-\sqrt{\frac{4+y}{2(y-2)}}\right)\right]+\frac{y+2}{2y}\\
&&-\frac{y+4}{2}\left\{\ln\left[4\left(1-\sqrt{\frac{4+y}{2(y-2)}}\right)+y\right]\right\}'+\left\{\ln\left[{4(y+4)\left(1-\sqrt{\frac{4+y}{2(y-2)}}\right)+2y}\right]\right\}'.
\end{eqnarray*}
Define
$$
r(y):=\frac{4}{y}\left(1-\sqrt{\frac{4+y}{2(y-2)}}\right).
$$
We have
\begin{eqnarray*}
-R'_4(y)
&\ge&\frac{1}{2}\left(r(y)-\frac{r^2(y)}{2}\right)-\frac{y+2}{2y}\\
&&+\frac{y+4}{2}\left\{\ln\left[4\left(1-\sqrt{\frac{4+y}{2(y-2)}}\right)+y\right]\right\}'-\left\{\ln\left[{4(y+4)\left(1-\sqrt{\frac{4+y}{2(y-2)}}\right)+2y}\right]\right\}'\\
&:=&Q_4(y).
\end{eqnarray*}
By virtue of Mathematica, we obtain that for $y\ge 15$,
\begin{eqnarray*}
&&y^2(4+y) \left(4+y-2\sqrt{\frac{2(4+y)}{y-2}}\right)\left[8-4\sqrt{\frac{2(4+y)}{y-2}}+y\left(3-\sqrt{\frac{2(4+y)}{y-2}}\right)\right]Q_4(y)\\
&=&4\left(\frac{4+y}{y-2}\right)^{\frac{3}{2}}\left[128\sqrt{\frac{4+y}{y-2}}+y\left(88\sqrt{2}-40\sqrt{\frac{4+y}{y-2}}\right)+y^2\left(26\sqrt{2}-34\sqrt{\frac{4+y}{y-2}}\right)\right.\\
&&\ \ \ \ \ \ \ \ \ \ \ \ \ \ \ \ \ \left.+ y^3\left(\sqrt{2}-\sqrt{\frac{4+y}{y-2}}\right)\right]\\
&>&4\left(\frac{4+y}{y-2}\right)^{\frac{3}{2}}\min\Bigg\{\min_{y\ge 17}\left[y^2\left(43\sqrt{2}-51\sqrt{\frac{4+y}{y-2}}\right)\right],\\
&&128\sqrt{\frac{4+15}{15-2}}+15\left(88\sqrt{2}-40\sqrt{\frac{4+15}{15-2}}\right)+15^2\left(26\sqrt{2}-34\sqrt{\frac{4+15}{15-2}}\right)+ 15^3\left(\sqrt{2}-\sqrt{\frac{4+15}{15-2}}\right)\Bigg\}
\end{eqnarray*}
\begin{eqnarray*}
&>&\frac{4}{51}\left(\frac{4+y}{y-2}\right)^{\frac{3}{2}}\min\Bigg\{\min_{y\ge 17}\left[y^2\left(1.19-\sqrt{1+\frac{6}{y-2}}\right)\right],19.87\Bigg\}\\
&\ge&\frac{4}{51}\left(\frac{4+y}{y-2}\right)^{\frac{3}{2}}\min\Bigg\{17^2\left(1.19-\sqrt{1+\frac{6}{17-2}}\right),19.87\Bigg\}\\
&>&1.96\cdot\frac{4}{51}\left(\frac{4+y}{y-2}\right)^{\frac{3}{2}}.
\end{eqnarray*}
Then,  $R'_4(y)<0$ for $15\le y\in\mathbb{R}$. Therefore, (\ref{RTD2}) holds for the case  $d_1=4$.

\subsection{Case $d_1=3$}

\subsubsection{Proof of  (\ref{May14z1})}

By (\ref{L1}), we get
\begin{eqnarray}\label{RTD3}
(\ref{May14z1})\ {\rm holds}
&\Leftarrow&2A(1-A)<d_2(B-A)\int_0^1\left(1-\frac{B-A}{1-A}x\right)^{\frac{d_2}{2}-1}dx\nonumber\\
&\Leftrightarrow&2A(1-A)<2(1-A)\left[1-\left(\frac{1-B}{1-A}\right)^{\frac{d_2}{2}}\right]\nonumber\\
&\Leftrightarrow&(1-B)^{\frac{d_2}{2}}<(1-A)^{\frac{d_2}{2}+1}.
\end{eqnarray}
Let $H_x(y)$ be defined by  (\ref{L12}). Note that  (\ref{L11}) and  (\ref{L13}) hold. Now let $x=3$. We have
\begin{eqnarray*}
H_3(y)=\left({\frac{y}{2}+1}\right)\ln\left[1+\frac{3}{y}\left(1+\sqrt{\frac{2(3+y)}{3(y-2)}}\right)\right].
\end{eqnarray*}
We will prove  (\ref{RTD3}) by showing that $H_3$ is a decreasing function for $y\in\{3,4,\dots,12\}$ and
$H'_3(y)<0$ for $12\le y\in\mathbb{R}$. By virtue of Mathematica, we obtain  the following values, which are decreasing.
\begin{eqnarray*}
&&H_3(3)=3.46574,\ \ H_3(4)=3.18962,\ \ H_3(5)=3.06414,\ \ H_3(6)=2.99125,\\
&&H_3(7)=2.94353,\ \ \ H_3(8)=2.9099,\ \ H_3(9)=2.88497,\ \ H_3(10)=2.86578,\\
&&H_3(11)=2.85058,\ \ H_3(12)=2.83826.
\end{eqnarray*}

We have
\begin{eqnarray*}
H'_3(y)=\frac{1}{2}\ln\left[1+\frac{3}{y}\left(1+\sqrt{\frac{2(3+y)}{3(y-2)}}\right)\right]+\left({\frac{y}{2}+1}\right)\left\{\ln\left[1+\frac{3}{y}\left(1+\sqrt{\frac{2(3+y)}{3(y-2)}}\right)\right]\right\}'.
\end{eqnarray*}
Define
$$
p_3(y):=\frac{3}{y}\left(1+\sqrt{\frac{2(3+y)}{3(y-2)}}\right).
$$
We have
\begin{eqnarray*}
H'_3(y)
&\le&\frac{1}{2}\left(p_3(y)-\frac{p_3^2(y)}{2}+\frac{p_3^3(y)}{3}\right)+\left({\frac{y}{2}+1}\right)\left\{\ln\left[1+\frac{3}{y}\left(1+\sqrt{\frac{2(3+y)}{3(y-2)}}\right)\right]\right\}'\\
&:=&L_3(y).
\end{eqnarray*}
By virtue of Mathematica and the fact
$$
\sqrt{6}>1.5\sqrt{\frac{3+y}{y-2}},\ \ \ \ y\ge 12,
$$
we obtain that for $y\ge 12$,
\begin{eqnarray*}
&&-4y^3(y-2)^2 \sqrt{\frac{3+y}{y-2}}\left(3+y+\sqrt{\frac{6(3+y)}{y-2}}\right)L_2(y)\\
&=&864\sqrt{\frac{3+y}{y-2}}-36y\left(11\sqrt{6}+6\sqrt{\frac{3+y}{y-2}}\right)-3y^2\left(29\sqrt{6}+88\sqrt{\frac{3+y}{y-2}}\right)\\
&&+y^3\left(19\sqrt{6}+9\sqrt{\frac{3+y}{y-2}}\right)+3y^4\left(\sqrt{6}-\sqrt{\frac{3+y}{y-2}}\right)\\
&>&y\left[(-396-87y+19y^2)\sqrt{6}+\left(-216-264y+9y^2+1.5y^3\right)\sqrt{\frac{3+y}{y-2}}\right].
\end{eqnarray*}

Define
\begin{eqnarray*}
U_{1}(y)&:=&-396-87y+19y^2,\\
U_{2}(y)&:=&2(-216-264y+9y^2+1.5y^3).
\end{eqnarray*}
We have
\begin{eqnarray*}
U_{1}(12+r)&=&1296 + 369 r + 19 r^2,\\
U_{2}(12+r)&=&1008 + 1200 r + 126 r^2 + 3 r^3.
\end{eqnarray*}
Then, $U_{1}(y),U_{2}(y)>0$ for $12\le y\in\mathbb{R}$. Thus, $H'_3(y)>0$ for $12\le y\in\mathbb{R}$. Therefore, (\ref{RTD3}) holds for the case $d_1=3$.

\subsubsection{Proof of  (\ref{May14z2})}

We assume without loss of generality that $D>C>0$. By (\ref{cond3}), we get $D>0\Leftrightarrow d_2\ge 15$. Further, by (\ref{78}), we get
\begin{eqnarray*}
D>C>0&\Leftrightarrow&3(d_2-2)(d_2+3)(d_2-4)^2>2d^2_2(d_2+1)^2\\
&\Leftrightarrow&-288 + 192 d_2 + 4 d_2^2 - 25 d_2^3 + d_2^4>0.
\end{eqnarray*}
Define
$$
P_3(d_2):=-288 + 192 d_2 + 4 d_2^2 - 25 d_2^3 + d_2^4.
$$
We have
\begin{eqnarray*}
&&P_3(15)=-30258,\ \ P_3(16)=-33056,\ \ P_3(17)=-35172,\ \ P_3(18)=-36360,\\
&&P_3(19)=-36350,\ \ P_3(20)=-34848,\ \ P_3(21)=-31536,\ \ P_3(22)=-26072,\\
&&P_3(23)=-18090,\ \ P_3(24)=-7200,\\
&&P(25+r)=7012 + 16017 r + 1879 r^2 + 75 r^3 + r^4.
\end{eqnarray*}
Then,
$$
D>C>0\Leftrightarrow d_2\ge 25.
$$

Similar to (\ref{L1}), we get
\begin{eqnarray}\label{May16hg}
&&(\ref{May14z2})\ {\rm holds}\nonumber\\
&\Leftrightarrow&2C(1-C)>d_2(D-C)\int_0^1\left(1+\frac{D-C}{C}x\right)^{0.5}\left(1-\frac{D-C}{1-C}x\right)^{\frac{d_2}{2}-1}dx\nonumber\\
&\Leftarrow&2C(1-C)>d_2(D-C)\int_0^1\left(1+\frac{D-C}{2C}x\right)\left(1-\frac{D-C}{1-C}x\right)^{\frac{d_2}{2}-1}dx\nonumber\\
&\Leftrightarrow&2C(1-C)>d_2(D-C)\left\{\frac{1+C}{2C}\int_0^1\left(1-\frac{D-C}{1-C}x\right)^{\frac{d_2}{2}-1}dx-\frac{1-C}{2C}\int_0^1\left(1-\frac{D-C}{1-C}x\right)^{\frac{d_2}{2}}dx\right\}\nonumber\\
&\Leftrightarrow&2C(1-C)>\frac{1-C^2}{C}\left[1-\left(\frac{1-D}{1-C}\right)^{\frac{d_2}{2}}\right]-\frac{d_2(1-C)^2}{(d_2+2)C}\left[1-\left(\frac{1-D}{1-C}\right)^{\frac{d_2}{2}+1}\right]\nonumber\\
&\Leftrightarrow&2C^2>\frac{2(d_2+1)C+2}{d_2+2}-\frac{2(1+C)+d_2(C+D)}{d_2+2}\left(\frac{1-D}{1-C}\right)^{\frac{d_2}{2}}\nonumber\\
&\Leftrightarrow&\left[2(1+C)+d_2(C+D)\right](1-D)^{\frac{d_2}{2}}>2[(d_2+2)C+1](1-C)^{\frac{d_2}{2}+1}\nonumber\\
&\Leftrightarrow&\ln\left[2(1+C)+d_2(C+D)\right]+{\frac{d_2}{2}}\ln(1-D)>\ln\left\{2[(d_2+2)C+1]\right\}+\left({\frac{d_2}{2}+1}\right)\ln(1-C)\nonumber\\
&\Leftrightarrow&\ln\left(1+\frac{C}{1-C}\right)>\ln\left(1+\frac{(d_2+2)C-d_2D}{2(1+C)+d_2(C+D)}\right)+\frac{d_2}{2}\ln\left(1+\frac{D-C}{1-D}\right).
\end{eqnarray}
Note that
\begin{eqnarray*}
&&(d_2+2)C>d_2D\\
&\Leftrightarrow&{\frac{d_2+2}{d_1+d_2\left(1-\sqrt{\frac{2(d_1+d_2)}{d_1(d_2-2)}}\right)^{-1}}}>{\frac{d_2}{d_1+(d_2-2)\left(1-\sqrt{\frac{2(d_1+d_2-2)}{d_1(d_2-4)}}\right)^{-1}}}\\
&\Leftrightarrow&\ \ 2d_1\left(1-\sqrt{\frac{2(d_1+d_2)}{d_1(d_2-2)}}\right)\left(1-\sqrt{\frac{2(d_1+d_2-2)}{d_1(d_2-4)}}\right)+(d_2^2-4)\left(1-\sqrt{\frac{2(d_1+d_2)}{d_1(d_2-2)}}\right)\\
&&>d_2^2\left(1-\sqrt{\frac{2(d_1+d_2-2)}{d_1(d_2-4)}}\right)\\
&\Leftrightarrow&\left(d_2^2-6+6\sqrt{\frac{2(3+d_2)}{3(d_2-2)}}\right)\sqrt{\frac{2(d_2+1)}{3(d_2-4)}}+2>(d_2^2+2)\sqrt{\frac{2(3+d_2)}{3(d_2-2)}}\\
&\Leftarrow&\left(d_2^2-6+2\sqrt{6}\right)\sqrt{\frac{2(d_2+1)}{3(d_2-4)}}>(d_2^2+2)\sqrt{\frac{2(3+d_2)}{3(d_2-2)}}\\
&\Leftarrow&(d_2^2-2)\sqrt{\frac{2(d_2+1)}{3(d_2-4)}}>(d_2^2+2)\sqrt{\frac{2(3+d_2)}{3(d_2-2)}}\\
&\Leftrightarrow&(d_2^2-2)^2(d_2+1)(d_2-2)>(d_2^2+2)^2(3+d_2)(d_2-4)\\
&\Leftrightarrow&2 (20 + 28 d_2^2 + 4 d_2^3 + d_2^4)>0.
\end{eqnarray*}
Then,  $(d_2+2)C>d_2D$. Note that $\ln(1+x)>x-\frac{x^2}{2}$ and $\ln(1+x)<x$ for any $x>0$. Thus, we get
\begin{eqnarray}\label{May17h}
&&\ln\left(1+\frac{C}{1-C}\right)>\ln\left(1+\frac{(d_2+2)C-d_2D}{2(1+C)+d_2(C+D)}\right)+\frac{d_2}{2}\ln\left(1+\frac{D-C}{1-D}\right)\nonumber\\
&\Leftarrow&\frac{C}{1-C}\left(1-\frac{C}{2(1-C)}\right)>\frac{(d_2+2)C-d_2D}{2(1+C)+d_2(C+D)}+\frac{d_2(D-C)}{2(1-D)}.
\end{eqnarray}

For $5\le y\in\mathbb{R}$, define
$$
C(y):={\frac{3}{3+y\left(1-\sqrt{\frac{2(3+y)}{3(y-2)}}\right)^{-1}}},\ \ \ \ D(y):={\frac{3}{3+(y-2)\left(1-\sqrt{\frac{2(y+1)}{3(y-4)}}\right)^{-1}}},
$$
$$
V(y):=\frac{C(y)}{1-C(y)}\left(1-\frac{C(y)}{2(1-C(y))}\right)-\frac{(y+2)C(y)-yD(y)}{2(1+C(y))+y(C(y)+D(y))}-\frac{y(D(y)-C(y))}{2(1-D(y))},
$$
and
$$
s_1(y):=\sqrt{\frac{y+1}{y-4}},\ \ \ \ s_2(y):=\sqrt{\frac{y+3}{y-2}}.
$$
By virtue of Mathematica, we get
$$
V(y)=\frac{G_1(y)}{G_2(y)},
$$
where
\begin{eqnarray*}
G_1(y)&=&-3\Bigg\{y^8\left[2\sqrt{6}s_1(y)-2\sqrt{6}s_2(y)\right] \\
&&\ \ \ \ \ \ +y^7\left [-4 - 8 \sqrt{6}s_1(y) +
   8 \sqrt{6}s_2(y) + 4s_1(y)s_2(y)\right]\\
&&\ \ \ \ \ \ +y^6 \left[-16 - 25  \sqrt{6}s_1(y) +
    23 \sqrt{6}s_2(y) -
    20  s_1(y)s_2(y)\right]\\
&&\ \ \ \ \ \ +y^5 \left[-392 + 77 \sqrt{6}s_1(y) +
 125 \sqrt{6}s_2(y)-
 70 s_1(y)s_2(y)\right]\\
&&\ \ \ \ \ \ +y^4 \left[-500 + 410 \sqrt{6}s_1(y) -
 364 \sqrt{6}s_2(y) -
 304 s_1(y)s_2(y)\right]\\
&&\ \ \ \ \ \ +4y^3 \left[1134 - 108 \sqrt{6}s_1(y)-
 457 \sqrt{6}s_2(y) +
 654 s_1(y)s_2(y)\right]\\
&&\ \ \ \ \ \ +16y^2 \left[435 - 237 \sqrt{6}s_1(y)+
 160\sqrt{6}s_2(y) +
 192s_1(y)s_2(y)\right]\\
&&\ \ \ \ \ \ +192y \left[-54 + 11 \sqrt{6}s_1(y) -
 66 s_1(y)s_2(y)\right]\\
&&\ \ \ \ \ \ +6912 \left[-1 + \sqrt{6}s_1(y)\right]\Bigg\},
\end{eqnarray*}
and
\begin{eqnarray*}
G_2(y)&=&2 (-4 + y)(-2 + y)^2) (y^2)  \left[3 + y -
  \sqrt{6}s_2(y)\right] \\
&&\cdot\Bigg\{y^2\left [-8 + \sqrt{6}s_1(y) + \sqrt{6}s_2(y)\right]+4\left[-3 + 3 \sqrt{6}s_1(y) +
\sqrt{6}s_2(y) -
 6  s_1(y)s_2(y)\right]\\
&&\ \ +2y \left[-13 + 4 \sqrt{6}s_1(y)+
 4 \sqrt{6}s_2(y)-
 6  s_1(y)s_2(y)\right]\Bigg\}.
\end{eqnarray*}

We will show that both $G_1(y)$ and $G_2(y)$ are negative for $y\in\{25,26,\dots\}$. Then $V(y)$ is a positive function of $y\in\{25,26,\dots\}$ and hence the proof will be complete by (\ref{May16hg}) and (\ref{May17h}).  Note that
$$
1<s_2(y)<s_1(y)\le\sqrt{\frac{26}{21}}\approx 1.112697,\ \ \ \ 25\le y\in\mathbb{R}.
$$
By virtue of Mathematica, we get
\begin{eqnarray*}
&&G_1(25)=-1.46179\cdot10^9,\ \ G_1(26)=-1.91825\cdot10^9,\ \ G_1(27)=-2.48438\cdot10^9,\\
&&G_1(28)=-3.17992\cdot10^9,\ \ G_1(29)=-4.02709\cdot10^9,\ \ G_1(30)=-5.05084\cdot10^9,\\
&&G_1(31)=-6.27904\cdot10^9,\ \ G_1(32)=-7.74269\cdot10^9,\ \ G_1(33)=-9.47617\cdot10^9.
\end{eqnarray*}
For  $25\le y\in\mathbb{R}$, we have
\begin{eqnarray*}
6912 \left[-1 + \sqrt{6}s_1(y)\right]>0,
\end{eqnarray*}
\begin{eqnarray*}
&&\ \ \ y^5 \left[-392 + 77 \sqrt{6}s_1(y) +
 125 \sqrt{6}s_2(y)-
 70 s_1(y)s_2(y)\right]\\
&&\ \ \ +192y \left[-54 + 11 \sqrt{6}s_1(y) -
 66 s_1(y)s_2(y)\right]\\
&&> y\Bigg\{25^4 \left[-392+202\sqrt{6}-70(1.113)^2\right]+192\left[-54 + 11 \sqrt{6} -
 66(1.113)^2\right]\Bigg\}\\
&&\approx (6.26157\cdot10^6)y>0,
\end{eqnarray*}
\begin{eqnarray*}
&&\ \ \ 16y^2 \left[435 - 237 \sqrt{6}s_1(y)+
 160\sqrt{6}s_2(y) +
 192s_1(y)s_2(y)\right]\\
&&>16y^2 \left[435 - 237 \sqrt{6}(1.113)+160\sqrt{6}+192\right]\approx  (16y^2)(372.7895)>0,
\end{eqnarray*}
and
\begin{eqnarray*}
&&\ \ \ 4y^3 \left[1134 - 108 \sqrt{6}s_1(y)-
 457 \sqrt{6}s_2(y) +
 654 s_1(y)s_2(y)\right]\\
&&>4y^3 \left[1134 -565\sqrt{6}(1.113)+654\right]\approx  (4y^3)(247.6506)>0.
\end{eqnarray*}
Note that
\begin{eqnarray*}
(y-5)s_1(y)s_2(y)-y>(y-5)[s_2(y)]^2-y=-\frac{15}{y-2},
\end{eqnarray*}
and, by setting $y=30+r$, we get
\begin{eqnarray*}
&&\sqrt{y-4}(2y^2-8y-23)>1.93\sqrt{y+1}( y^2- 6 y +8) \\
&\Leftrightarrow&-23543936 - 8238032 y + 6438956 y^2 - 489960 y^3 - 70261 y^4 +
 2751 y^5>0\\
&\Leftrightarrow&2233345504 + 2608569328 r + 325703156 r^2 + 15837720 r^3 +
 342389 r^4 + 2751 r^5>0,
\end{eqnarray*}
which holds for $r\ge 0$. Then, we obtain that for $34\le y\in\mathbb{R}$,
\begin{eqnarray*}
&&\ \ \ y^8\left[2\sqrt{6}s_1(y)-2\sqrt{6}s_2(y)\right]+y^7\left [-4 - 8 \sqrt{6}s_1(y) +
   8 \sqrt{6}s_2(y) + 4s_1(y)s_2(y)\right]\\
&&\ \ \ +y^6 \left[-16 - 25  \sqrt{6}s_1(y) +
    23 \sqrt{6}s_2(y) -
    20  s_1(y)s_2(y)\right] \\
&&\ \ \ +y^4 \left[-500 + 410 \sqrt{6}s_1(y) -
 364 \sqrt{6}s_2(y) -
 304 s_1(y)s_2(y)\right]\\
&&=\sqrt{6}y^6\left[s_1(y)-s_2(y)\right](2y^2-8y-23)-y^6\left\{16+2\sqrt{6}s_1(y)-4\left[(y-5)s_1(y)s_2(y)-y\right]\right\}\\
&&\ \ \ +y^4 \left[-500 + 410 \sqrt{6}s_1(y) -
 364 \sqrt{6}s_2(y) -
 304 s_1(y)s_2(y)\right]\\
&&>\frac{10\sqrt{6}y^6(2y^2-8y-23)}{[s_1(y)+s_2(y)]( y^2- 6 y +8) }-y^6\left[16+2\sqrt{6}(1.113)+\frac{60}{y-2}\right]\\
&&\ \ \ -y^4\left[500-410 \sqrt{6}+364(1.113)+304(1.113)^2\right]\\
&&>\frac{5\sqrt{6}y^6\sqrt{y-4}(2y^2-8y-23)}{\sqrt{y+1}( y^2- 6 y +8) }-y^6\left[16+2\sqrt{6}(1.113)+\frac{60}{y-2}\right]-278y^4\\
&&>y^6\left[5\sqrt{6}\cdot 1.93-16-2\sqrt{6}(1.113)-\frac{60}{y-2}-\frac{278}{y^2}\right]\\
&&\ge y^6\left[5\sqrt{6}\cdot1.93-16-2\sqrt{6}(1.113)-\frac{60}{34-2}-\frac{278}{34^2}\right]\\
&&\approx 0.0695274
\end{eqnarray*}
Thus, $G_1(y)<0$ for $y\in\{25,26,\dots\}$.

For $25\le y\in\mathbb{R}$, define
$$
y=25+x,\ \ \ \ u=\sqrt{\frac{26+x}{21+x}},\ \ \ \ v=\sqrt{\frac{28+x}{23+x}}.
$$
We have that
$$
1<v<u\le\sqrt{\frac{26}{21}}\approx 1.112697.
$$
By virtue of Mathematica, we get
\begin{eqnarray*}
G_2(y)&=&2 (25 + x)^2 (483 + 44 x + x^2) \\
&&\cdot\Bigg\{-3785600 + 548100 \sqrt{6}u +
 664102 \sqrt{6}v-
 324162 uv+x^4 \left[-8 + \sqrt{6}u + \sqrt{6}v\right]\\
&&\ \  +x (-577824 + 80699  \sqrt{6}u+
   95091  \sqrt{6}v -
   37278  uv) \\
&&\ \ +x^2 (-33056 + 4451  \sqrt{6}u +
   5041  \sqrt{6}v-
   1422  uv)\\
&&\ \ +x^3 (-840 + 109 \sqrt{6}u+
   117\sqrt{6}v-
   18 uv)\Bigg\}\\
&<&2 (25 + x)^2 (483 + 44 x + x^2) \\
&&\cdot\Bigg\{-3785600 + 548100 \sqrt{6}(1.113) +
 664102 \sqrt{6}(1.113)-
 324162\\
&&\ \  +x^4 \left[-8 + \sqrt{6}(1.113) + \sqrt{6}(1.113)\right]\\
&&\ \  +x (-577824 + 80699  \sqrt{6}(1.113)+
   95091  \sqrt{6}(1.113) -
   37278  ) \\
&&\ \ +x^2 (-33056 + 4451  \sqrt{6}(1.113) +
   5041  \sqrt{6}(1.113)-
   1422  )\\
&&\ \ +x^3 (-840 + 109 \sqrt{6}(1.113)+
   117\sqrt{6}(1.113)-
   18 )\Bigg\}\\
&<&0.
\end{eqnarray*}

\section{Remarks on proof of  (\ref{30A}) for cases $d_1\ge 5$}\setcounter{equation}{0}

\subsection{Case $d_1=2(k+1)$, $2\le k\in\mathbb{N}$}\setcounter{equation}{0}

In this subsection, we assume that $d_1\ge 6$ and $d_1$ is an even number. By (\ref{May1311}), it suffices to prove inequalities (\ref{May14z1}) and (\ref{May14z2}).

By (\ref{May1311}), we get
\begin{eqnarray}\label{May15d1}
&&(\ref{May14z1})\ {\rm holds}\nonumber\\
&\Leftrightarrow&I_A\left(\frac{d_1}{2},\frac{d_2+2}{2}\right)<I_B\left(\frac{d_1}{2},\frac{d_2}{2}\right)\nonumber\\
&\Leftrightarrow&\frac{\int_0^At^{\frac{d_1}{2}-1}(1-t)^{\frac{d_2+2}{2}-1}dt}{B\left(\frac{d_1}{2},\frac{d_2+2}{2}\right)}<\frac{\int_0^Bt^{\frac{d_1}{2}-1}(1-t)^{\frac{d_2}{2}-1}dt}{B\left(\frac{d_1}{2},\frac{d_2}{2}\right)}\nonumber\\
&\Leftrightarrow&\frac{\int_{1-A}^1(1-t)^{\frac{d_1}{2}-1}t^{\frac{d_2+2}{2}-1}dt}{B\left(\frac{d_1}{2},\frac{d_2+2}{2}\right)}<\frac{\int_{1-B}^1(1-t)^{\frac{d_1}{2}-1}t^{\frac{d_2}{2}-1}dt}{B\left(\frac{d_1}{2},\frac{d_2}{2}\right)}\nonumber\\
&\Leftrightarrow&\frac{\int_0^{1-A}(1-t)^{\frac{d_1}{2}-1}t^{\frac{d_2+2}{2}-1}dt}{B\left(\frac{d_1}{2},\frac{d_2+2}{2}\right)}>\frac{\int_0^{1-B}(1-t)^{\frac{d_1}{2}-1}t^{\frac{d_2}{2}-1}dt}{B\left(\frac{d_1}{2},\frac{d_2}{2}\right)}\nonumber\\
&\Leftrightarrow&(d_1+d_2)\int_0^{1-A}(1-t)^{\frac{d_1}{2}-1}t^{\frac{d_2+2}{2}-1}dt>d_2\int_0^{1-B}(1-t)^{\frac{d_1}{2}-1}t^{\frac{d_2}{2}-1}dt\nonumber\\
&\Leftrightarrow&(d_1+d_2)\int_0^{1-A}\sum_{n=0}^{\frac{d_1}{2}-1}(-1)^{n}{\frac{d_1}{2}-1\choose n}t^{n+\frac{d_2+2}{2}-1}dt>d_2\int_0^{1-B}\sum_{n=0}^{\frac{d_1}{2}-1}(-1)^{n}{\frac{d_1}{2}-1\choose n}t^{n+\frac{d_2}{2}-1}dt\nonumber\\
&\Leftrightarrow&(d_1+d_2)(1-A)^{\frac{d_2+2}{2}}\sum_{n=0}^{\frac{d_1}{2}-1}{\frac{d_1}{2}-1\choose n}\frac{(A-1)^{n}}{2n+d_2+2}>d_2(1-B)^{\frac{d_2}{2}}\sum_{n=0}^{\frac{d_1}{2}-1}{\frac{d_1}{2}-1\choose n}\frac{(B-1)^{n}}{2n+d_2}\nonumber\\
&\Leftrightarrow&\ \ (1-B)^{\frac{d_2}{2}}\left\{\prod_{n=0}^{\frac{d_1}{2}-1}(2n+d_2)\right\}\sum_{n=0}^{\frac{d_1}{2}-1}{\frac{d_1}{2}-1\choose n}\frac{(B-1)^{n}}{2n+d_2}\nonumber\\
&&<(1-A)^{\frac{d_2+2}{2}}\left\{\prod_{n=0}^{\frac{d_1}{2}-1}(2n+d_2+2)\right\}\sum_{n=0}^{\frac{d_1}{2}-1}{\frac{d_1}{2}-1\choose n}\frac{(A-1)^{n}}{2n+d_2+2}.
\end{eqnarray}
For $d_1\in\mathbb{N}$ and $3\le y\in\mathbb{R}$, define
\begin{eqnarray*}
J_{d_1}(y)&:=&-\frac{y+2}{2}\ln\left[1+\frac{d_1}{y}\left(1+\sqrt{\frac{2(d_1+y)}{d_1(y-2)}}\right)\right]+\sum_{n=0}^{\frac{d_1}{2}-1}\ln(2n+y+2)\\
&&+\ln\left\{\sum_{n=0}^{\frac{d_1}{2}-1}{\frac{d_1}{2}-1\choose n}\frac{(-1)^{n}}{(2n+y+2)\left[1+\frac{d_1}{y}\left(1+\sqrt{\frac{2(d_1+y)}{d_1(y-2)}}\right)\right]^n}\right\}.
\end{eqnarray*}
We may try to show that $J_{d_1}$ is a strictly increasing function by proving that
$J'_{d_1}(y)>0$.

We assume without loss of generality that $D>C>0$. By (\ref{78}), this assumption is equivalent to
\begin{eqnarray*}
d_1(d_2-2)(d_1+d_2)(d_2-4)^2>2d^2_2(d_1+d_2-2)^2.
\end{eqnarray*}
Similar to (\ref{May15d1}), we can show that
\begin{eqnarray*}
&&(\ref{May14z2})\ {\rm holds}\nonumber\\
&\Leftrightarrow&\ \ (1-D)^{\frac{d_2}{2}}\left\{\prod_{n=0}^{\frac{d_1}{2}-1}(2n+d_2)\right\}\sum_{n=0}^{\frac{d_1}{2}-1}{\frac{d_1}{2}-1\choose n}\frac{(D-1)^{n}}{2n+d_2}\nonumber\\
&&>(1-C)^{\frac{d_2+2}{2}}\left\{\prod_{n=0}^{\frac{d_1}{2}-1}(2n+d_2+2)\right\}\sum_{n=0}^{\frac{d_1}{2}-1}{\frac{d_1}{2}-1\choose n}\frac{(C-1)^{n}}{2n+d_2+2}.
\end{eqnarray*}
For $d_1\in\mathbb{N}$ and $3\le y\in\mathbb{R}$, define
\begin{eqnarray*}
K_{d_1}(y)&:=&-\frac{y+2}{2}\ln\left[1+\frac{d_1}{y}\left(1-\sqrt{\frac{2(d_1+y)}{d_1(y-2)}}\right)\right]+\sum_{n=0}^{\frac{d_1}{2}-1}\ln(2n+y+2)\\
&&+\ln\left\{\sum_{n=0}^{\frac{d_1}{2}-1}{\frac{d_1}{2}-1\choose n}\frac{(-1)^{n}}{(2n+y+2)\left[1+\frac{d_1}{y}\left(1-\sqrt{\frac{2(d_1+y)}{d_1(y-2)}}\right)\right]^n}\right\}.
\end{eqnarray*}
We may try to show that $K_{d_1}$ is a strictly decreasing function by proving that
$K'_{d_1}(y)<0$.

\subsection{Case $d_1=2k+1$, $2\le k\in\mathbb{N}$}

In this subsection, we assume that $d_1\ge 5$ and $d_1$  is an odd number. By (\ref{May1311}), it suffices to prove (\ref{May14z1}) and (\ref{May14z2}).

For $x>0$, denote by  $\lfloor x\rfloor$  the greatest integer less than or equal to $x$ and
$$
(x)_{\langle n\rangle}:=x(x-1)\cdots(x-n+1).
$$
Similar to (\ref{L1}), we get
\begin{eqnarray}\label{RT}
&&(\ref{May14z1})\ {\rm holds}\nonumber\\
&\Leftrightarrow&2A(1-A)<d_2(B-A)\int_0^1\left(1+\frac{B-A}{A}x\right)^{\frac{d_1}{2}-1}\left(1-\frac{B-A}{1-A}x\right)^{\frac{d_2}{2}-1}dx\nonumber\\
&\Leftarrow&2A(1-A)<d_2(B-A)\sum_{n=0}^{\lfloor \frac{d_1}{2}-1\rfloor}\frac{\left(\frac{d_1}{2}-1\right)_{\langle n\rangle}\left(\frac{B-A}{A}\right)^n}{n!}\int_0^1x^n\left(1-\frac{B-A}{1-A}x\right)^{\frac{d_2}{2}-1}dx\nonumber\\
&\Leftrightarrow&2A<d_2\sum_{n=0}^{\lfloor \frac{d_1}{2}-1\rfloor}\frac{\left(\frac{d_1}{2}-1\right)_{\langle n\rangle}\left(\frac{1-A}{A}\right)^n}{n!}\sum_{l=0}^n\frac{(-1)^l}{\frac{d_2}{2}+l}{n\choose l}\left[1-\left(\frac{1-B}{1-A}\right)^{\frac{d_2}{2}+l}\right].
\end{eqnarray}

We assume without loss of generality that $D>C>0$. By (\ref{78}), this assumption is equivalent to
\begin{eqnarray*}
d_1(d_2-2)(d_1+d_2)(d_2-4)^2>2d^2_2(d_1+d_2-2)^2.
\end{eqnarray*}
For $x>0$, denote by  $\lceil x\rceil$  the smallest integer greater than or equal to $x$. Similar to (\ref{RT}), we get
\begin{eqnarray}\label{ASD}
&&(\ref{May14z2})\ {\rm holds}\nonumber\\
&\Leftarrow&2C>d_2\sum_{n=0}^{\lceil  \frac{d_1}{2}-1\rceil}\frac{\left(\frac{d_1}{2}-1\right)_{\langle n\rangle}\left(\frac{1-C}{C}\right)^n}{n!}\sum_{l=0}^n\frac{(-1)^l}{\frac{d_2}{2}+l}{n\choose l}\left[1-\left(\frac{1-D}{1-C}\right)^{\frac{d_2}{2}+l}\right].
\end{eqnarray}
Conjecture \ref{con1} can be resolved if we are able to establish inequalities (\ref{RT}) and (\ref{ASD}) with two variables $d_1=2k+1$, $2\le k\in\mathbb{N}$, and $d_2\in\mathbb{N}$.

\vskip 0.8cm
{ \noindent {\bf\large Acknowledgements}\quad  This work was supported by the National Natural Science Foundation of China (No. 12171335), the Science Development Project of Sichuan University (No. 2020SCUNL201) and the Natural Sciences and Engineering Research Council of Canada (No. 4394-2018).

\end{document}